\documentclass[a4paper,12pt, twoside, reqno]{amsart}
\usepackage{amsmath}%
\usepackage{amsfonts}%
\usepackage{amssymb}%
\usepackage{graphicx}
\newtheorem{theorem}{Theorem}
\theoremstyle{plain}

\newtheorem{example}{Example}

\newtheorem{remark}{Remark}

\numberwithin{equation}{section}

\begin{document}
\title{Coupled fixed point theorems for generalized symmetric Meir–-Keeler contractions in ordered metric spaces}
\author{Vasile Berinde}

\begin{abstract}
In this paper we introduce generalized symmetric Meir–Keeler contractions  and prove some coupled fixed point theorems for mixed monotone operators $F:X \times X \rightarrow  X$ in partially ordered metric spaces. The obtained results extend, complement and unify some recent coupled fixed point theorems due to Samet [B. Samet, \textit{Coupled fixed point theorems for a generalized Meir–-Keeler contraction
in partially ordered metric spaces}, Nonlinear Anal.  \textbf{72} (2010), 4508--4517], Bhaskar and Lakshmikantham [T.G. Bhaskar, V. Lakshmikantham, \textit{Fixed point theorems in partially ordered metric spaces and applications}, Nonlinear Anal. TMA \textbf{65} (2006) 1379-1393] and some other very recent papers. An example to show that our generalizations are effective is also presented.
\end{abstract}
\maketitle

\pagestyle{myheadings} \markboth{Vasile Berinde} {Coupled fixed point theorems for generalized symmetric Meir–Keeler contractions}

\section{Introduction} 
Let $(X,d)$ be a metric space and $T:X\rightarrow X$ a self mapping.  If $(X,d)$ is complete and $T$ is a contraction, i.e., there exists a constant $\alpha \in [0,1)$ such that 
\begin{equation} \label{eq1.1}
d(Tx, Ty)\leq a\,d(x,y),\text{ for all}\;x,y\in X,
\end{equation}
then, by Banach contraction mapping principle, which is a classical and powerful tool in nonlinear analysis, we know that $T$ has a unique fixed point $p$ and, for any $x_0\in X$, the Picard iteration $\{T^n x_{0}\}$ converges to $p$.

The Banach contraction mapping principle has been generalized in several directions, see for example \cite{Ber07} and \cite{Rus2} for  recent surveys. 
One of these generalizations, known as the Meir-Keeler fixed point theorem \cite{Meir}, has been obtained by replacing the contraction condition \eqref{eq1.1} by the following more general assumption: for all $\epsilon>0$ there exists $\delta(\epsilon)>0$ such that
\begin{equation} \label{eq1.2}
x,y\in X,\,\epsilon\leq d(x,y) <\epsilon+\delta(\epsilon)\, \Rightarrow\, (Tx, Ty)<\epsilon.
\end{equation}

Recently, Ran  and  Reurings  \cite{Ran} initiated an other important direction in generalizing the Banach contraction mapping principle, by considering a partial ordering on the metric space $(X,d)$ and by requiring that the contraction condition \eqref{eq1.1} is satisfied only for comparable elements, that is, we have
\begin{equation} \label{eq1.3}
d(Tx, Ty)\leq a\,d(x,y),\textnormal{ for all}\;x,y\in X,\,\textnormal{ with } x \geq y.
\end{equation}
In compensation, the authors in \cite{Ran} assumed that $T$ satisfies a certain monotonicity condition. 

This new approach has been then followed by several authors: Agarwal  et  al.  \cite{Aga},  Nieto  and  Lopez  \cite{Nie06}, \cite{Nie07}, O'Regan and Petru\c sel \cite{Oreg},  who obtained fixed point theorems, and also by Bhaskar  and  Lakshmikantham  \cite{Bha}, Lakshmikantham and Ciric \cite{LakC}, Luong and Thuan \cite{Luong}, Samet \cite{Sam} and many others, who obtained  coupled fixed point theorems or coincidence point theorems. These results also found important applications to the existence of solutions for matrix equations or ordinary differential equations and integral equations, see  \cite{Bha}, \cite{Luong}, \cite{Nie06}, \cite{Nie07}, \cite{Oreg}, \cite{Ran}, and some of the references therein.

In order to state the main result in \cite{Sam}, we recall the following notions. Let $\left(X,\leq\right)$ be a partially ordered set and  endow the product space $X \times X$ with the following partial order:$$\textnormal{for } \left(x,y\right), \left(u,v\right) \in X \times X,  \left(u,v\right) \leq \left(x,y\right) \Leftrightarrow x\geq u,  y\leq v.$$

We say that a mapping $F:X \times X \rightarrow X$ has the \textit{mixed monotone property} if $F\left(x,y\right)$ is monotone nondecreasing in $ x$ and is monotone non increasing in $y$, that is, for any $x, y \in X,$   
\begin{equation} \label{eq1.4}
x_{1}, x_{2} \in X,  x_{1} \leq  x_{2} \Rightarrow  F\left(x_{1},y\right) \leq  F\left(x_{2},y\right) 
\end{equation}
  and, respectively,
\begin{equation} \label{eq1.5}
 y_{1},  y_{2} \in  X, y_{1} \leq y_{2} \Rightarrow F\left(x,y_{1}\right)  \geq F\left(x,y_{2}\right). 
\end{equation}

 We say F has the \textit{strict mixed monotone property} if the strict inequality in the left hand side of \eqref{eq1.4} and \eqref{eq1.5} implies the strict inequality in the right hand side, respectively.

A pair $ \left(x,y\right) \in X \times X $ is called a \textit{coupled  fixed  point} of  $F$ if                                              $$ F\left(x,y\right) = x \textnormal{ and } F\left(y,x\right) = y.$$
The next theorem is the main existence result in \cite{Sam}.

\begin{theorem}[Samet \cite{Sam}]\label{th1}
	Let  $\left(X,\leq\right)$ be a partially ordered set and suppose there is a metric $d$  on $X$ such that $\left(X,d\right)$  is a complete metric space. Let $F : X \times X \rightarrow X $ be a continuous mapping having the strict mixed monotone property on $X$. Assume also that $F$ is a generalized Meir-Keeler operator, that is, for each $\epsilon>0$, there exists $\delta(\epsilon)>0$ such that      
\begin{equation} \label{Samet}
\epsilon\leq \frac{1}{2}\left[d\left(x,u\right) + d\left(y,v\right)\right]<\epsilon+\delta(\epsilon)\,\Rightarrow\,	d\left(F\left(x,y\right),F\left(u,v\right)\right) <\epsilon, 
\end{equation} 
for all $x,y\in X$ satisfying $x \geq u, y \leq  v$.

If there exist $x_{0}, y_{0} \in X$ such that  
$$                                                                                       
x_{0} < F\left(x_{0},y_{0}\right)\textrm{ and }y_{0} \geq F\left(y_{0},x_{0}\right),
$$
then there exist $x, y \in X$ such that $$x = F\left(x,y\right)\textnormal{ and }y = F\left(y,x\right).$$ 
	\end{theorem}        
	
In the same paper \cite{Sam} the author also established other existence as well as existence and uniqueness results for coupled fixed points of mixed strict monotone generalized Meir-Keeler operators.

Starting from the results in \cite{Sam}, our main aim in this paper is to obtain more general coupled fixed point theorems for mixed monotone operators $F : X \times X \rightarrow X $  satisfying a generalized Meir-Keeler contractive condition which is significantly weaker that \eqref{Samet}. Our technique of proof is different and slightly simpler than the ones used in \cite{Sam} and \cite{Bha}. 
On the other hand, we do not require that $F$ is continuous. We thus extend, unify generalize and complement several related results in literature, amongst which we mention the ones in \cite{Aga}, \cite{Bha}, \cite{Meir}, \cite{Ran} and \cite{Sam}.

\section{Main results}
	
The first main result in this paper is the following coupled fixed point result which generalize Theorem \ref{th1} (Theorem 2.1 in \cite{Sam}) and Theorem 2.1 in \cite{Bha} and some other related results.

\begin{theorem}\label{th2}
	Let  $\left(X,\leq\right)$ be a partially ordered set and suppose there is a metric $d$  on $X$ such that $\left(X,d\right)$  is a complete metric space. Assume $F : X \times X \rightarrow X $ has the mixed monotone property and is also a generalized symmetric Meir-Keeler operator, that is, for each $\epsilon>0$, there exists $\delta(\epsilon)>0$ such that  for all $x,y\in X$  satisfying $x \geq u, y \leq  v$,
$$
\epsilon\leq \frac{1}{2}\left[d\left(x,u\right) + d\left(y,v\right)\right]<\epsilon+\delta(\epsilon)
$$     
implies                              
\begin{equation} \label{Samet1}
\frac{1}{2}\left[d\left(F\left(x,y\right),F\left(u,v\right)\right)+d\left(F\left(y,x\right),F\left(v,u\right)\right)\right] < \epsilon.
\end{equation} 
  
If there exist $x_{0}, y_{0} \in X$ such that  
\begin{equation} \label{mic}
 x_{0} \leq F\left(x_{0},y_{0}\right)\textrm{ and }y_{0} \leq F\left(y_{0},x_{0}\right),
\end{equation} or
\begin{equation} \label{mare} 
x_{0} \geq F\left(x_{0},y_{0}\right)\textrm{ and }y_{0} \leq F\left(y_{0},x_{0}\right),
\end{equation}
 then there exist $\overline{x}, \overline{y} \in X$ such that $$\overline{x} = F\left(\overline{x},\overline{y}\right)\textrm{and }\overline{y} = F\left(\overline{y},\overline{x}\right).$$ 
	\end{theorem}

\begin{proof}
Consider the functional $d_2:X^2\times X^2 \rightarrow \mathbb{R}_{+}$ defined by
$$
d_2(Y,V)=\frac{1}{2}\left[d(x,u)+d(y,v)\right],\,\forall Y=(x,y),V=(u,v) \in X^2.
$$
It is a simple task to check that $d_2$ is a metric on $X^2$ and, moreover, that, if $(X,d)$ is complete, then $(X^2,d_2)$ is a complete metric space, too.
Now consider the operator $T:X^2\rightarrow X^2$ defined by
$$
T(Y)=\left(F(x,y),F(y,x)\right),\,\forall Y=(x,y) \in X^2.
$$
Clearly, for $Y=(x,y),\,V=(u,v)\in X^2$, in view of the definition of $d_2$, we have
$$
d_2(T(Y),T(V))=\frac{d\left(F\left(x,y\right),F\left(u,v\right)\right)+d\left(F\left(y,x\right),F\left(v,u\right)\right)}{2}
$$
and
$$
d_2(Y,V)=\frac{d\left(x,u\right) + d\left(y,v\right)}{2}.
$$
Hence, by the contractive condition \eqref{Samet1} we obtain a usual Meir-Keeler type condition: for all $\epsilon>0$ there exists $\delta(\epsilon)>0$ such that
\begin{equation} \label{contr}
Y,V\in X^2,\,Y\geq V,\,\epsilon\leq d_2(Y,V) <\epsilon+\delta(\epsilon)\, \Rightarrow\, d_2(T(Y), T(V))<\epsilon.
\end{equation}
Assume \eqref{mic} holds (the case \eqref{mare} is similar). Then, there exists $x_0,y_0\in X$ such that
$$
x_0\leq F(x_0,y_0) \textnormal{ and } y_0\geq F(y_0,x_0). 
$$
Denote $Z_0=(x_0,y_0)\in X^2$ and consider the Picard iteration associated to $T$ and to the initial approximation $Z_0$, that is, the sequence $\{Z_n\}\subset X^2$ defined by
\begin{equation} \label{eq-3}
Z_{n+1}=T (Z_n),\,n\geq 0,
\end{equation}
where $Z_n=(x_n,y_n)\in X^2,\,n\geq 0$.

Since $F$ is mixed monotone, we have
$$
Z_0=(x_0,y_0)\leq (F(x_0,y_0), F(y_0,x_0))=(x_1,y_1)=Z_1
$$
and, by induction,
$$
Z_n=(x_n,y_n)\leq (F(x_n,y_n), F(y_n,x_n))=(x_{n+1},y_{n+1})=Z_{n+1},
$$
which shows that $T$ is monotone and the sequence $\{Z_n\}_{n=0}^{\infty}$ is non decreasing.

Note that \eqref{contr} implies the strict contractive condition
\begin{equation} \label{eq-5.5}
d_2(T(Y),T(Z)< d_2 (Y, Z),\, Y>Z,
\end{equation}
Take now $Y=Z_n > Z_{n-1}=V$ in \eqref{eq-5.5} to obtain
$$
d_2(T(Z_{n}),T(Z_{n-1})< d_2 (Z_n, Z_{n-1}),\,n\geq 1,
$$
which shows that the sequence of nonnegative numbers $\{\eta_n\}_{n=0}^{\infty}$ given by
\begin{equation} \label{eq-4}
\eta_n=d_2(Z_{n},Z_{n-1}),\,n\geq 1,
\end{equation}
is non-increasing, hence convergent to some $\epsilon\geq 0$.  

We now prove that necessarily $\epsilon=0$. Suppose, to the contrary, that $\epsilon>0$. Then, there exist a positive integer $p$ such that
$$
\epsilon\leq \eta_p<\epsilon+\delta(\epsilon),
$$
which, by the Meir-Keeler condition \eqref{contr}, yields
$$
\eta_{p+1}=d_2(T(Z_{p}),T(Z_{p-1}))<\epsilon,
$$
a contradiction, since $\{\eta_n\}_{n=0}^{\infty}$  converges non-increasingly to $\epsilon$. Therefore $\epsilon=0$, that is, 
\begin{equation} \label{eq-5}
\lim_{n\rightarrow \infty} \eta_n=\frac{1}{2}\lim_{n\rightarrow \infty}\left[d(x_{n+1},x_n)+d(y_{n+1},y_n)\right]=0.
\end{equation}
Let now $\epsilon>0$ be arbitrary and $\delta(\epsilon)$ the corresponding value from the hypothesis of our Theorem. By \eqref{eq-5}, there exists a positive integer $k$ such that 
\begin{equation} \label{eq-5.1}
\eta_{n+1}=d_2(Z_{k+1}, Z_k)<\delta(\epsilon).
\end{equation}
For this fixed number $k$, consider now the set
$$
\Lambda_k:=\{Z=(x,y)\in X^2: x>x_k,y\leq y_k \textnormal{ and } d_2(Z_k,Z)<\epsilon+\delta(\epsilon)\}.
$$
By \eqref{eq-5.1}, $\Lambda_k\neq \emptyset$. We claim that
\begin{equation} \label{eq-5.2}
Z\in \Lambda_k \Rightarrow T(Z)\in \Lambda_k.
\end{equation}  
Indeed, let $Z\in \Lambda_k$. Then $d_2(Z_k,Z)<\epsilon+\delta(\epsilon)$ and hence
$$
d_2(Z_k,T(Z))\leq d_2(Z_k,T(Z_k)))+d_2(T(Z_k),T(Z))=
$$
$$
=d_2(Z_k,Z_{k+1})+d_2(T(Z_k),T(Z))
$$
which, by \eqref{eq-5.2} and Meir-Keeler type condition \eqref{Samet1}, is
$$
<\delta(\epsilon)+\epsilon.
$$
Thus, by \eqref{eq-5.2} we have $Z_{k+1}\in \Lambda_k$ and, by induction,
$$
Z_{n}=(x_n,y_n)\in \Lambda_k, \, \textnormal{ for all } n>k.
$$
This implies that for all $n,m>k$, we have
$$
d_2(Z_n,Z_m)\leq d_2(Z_n,Z_k)+d_2(Z_m,Z_k)<2(\epsilon+\delta(\epsilon))\leq  4 \epsilon.
$$
Therefore $\{Z_n\}_{n=0}^{\infty}$ is a Cauchy sequence in the complete metric space $(X^2,d_2)$ and hence there exists $\overline{Z}\in X^2$ such that
$$
\lim_{n\rightarrow \infty} Z_n=\overline{Z}.
$$ 
 
By condition \eqref{eq-5.1}, $T$ is continuous in $(X^2,d_2)$, and hence by \eqref{eq-3} it follows that $\overline{Z}$ is a fixed point of $T$, that is,
$$
T(\overline{Z})=\overline{Z}.
$$
Let $\overline{Z}=(\overline{x},\overline{y})$. Then, by the definition of $T$, this means
$$\overline{x} = F\left(\overline{x},\overline{y}\right)\textrm{and }\overline{y} = F\left(\overline{y},\overline{x}\right),$$ 
that is,  $(\overline{x},\overline{y})$ is a coupled fixed point of $F$.
\end{proof}

 \begin{remark} \em
Theorem \ref{th2} is more general than Theorem \ref{th1} (i.e., Theorem 2.1 in \cite{Sam}), since the contractive condition \eqref{Samet1} is weaker than \eqref{Samet}, a fact which is clearly illustrated by Example \ref{ex1}.

Secondly, while Theorem 2.1 in \cite{Sam} assumed that $F$ is continuous, in our Theorem \ref{th2} we did not use this assumption. 

Apart from these improvements, we note that our proof is significantly simpler and shorter than the one in \cite{Sam}.
 \end{remark}
 
\begin{example} \label{ex1} \em
Let $X=\mathbb{R},$ $d\left(x,y\right)=|x-y|$ and $F:X\times X \rightarrow X$ be defined by 
$$
F\left(x,y\right)=\frac{x-3y}{5}, \,(x,y)\in X^2.
$$ 
Then $F$ is mixed monotone and satisfies condition \eqref{Samet1}  but  does not satisfy condition \eqref{Samet}. 

Assume, to the contrary, that  \eqref{Samet} holds. Let $x,y,u,v \in X$, $x\geq u,\,y\leq v$, such that
$$
\epsilon \leq \frac{1}{2}\left[\left|x-u\right|+\left|y-v\right|\right]<\epsilon+\delta(\epsilon).
$$
For $x=u$, this gives
\begin{equation} \label{eq-5.3}
\epsilon\leq  \frac{\left|y-v\right|}{2}<\epsilon+\delta(\epsilon),\,y\leq v.
\end{equation} 
which by \eqref{Samet} would imply
$$
\left|\frac{x-3y}{5}-\frac{u-3v}{5}\right|= \left|\frac{3v-3y}{5}\right|=\dfrac{3}{5}\left|y-v\right|<\epsilon,\,x= u,\,y< v,
$$
and this in turn, by \eqref{eq-5.3},  would imply
$$
2 \epsilon\leq \left|y-v\right|  <\dfrac{5}{3}\cdot \epsilon < 2\epsilon,
$$
a contradiction. Hence $F$  does not satisfy condition \eqref{Samet}.

Now we prove that \eqref{Samet1} holds. Indeed, we have
$$
\left|\frac{x-3y}{5}-\frac{u-3v}{5}\right| \leq \frac{1}{5}\left|x-u\right|+\frac{3}{5}\left|y-v\right|,\,x\geq u,\,y\leq v,
$$
and
$$
\left|\frac{y-3x}{5}-\frac{v-3u}{5}\right| \leq \frac{1}{5}\left|y-v\right|+\frac{3}{5}\left|x-u\right|,\,x\geq u,\,y\leq v,
$$
and by summing up the two inequalities above we get for all $\,x\geq u,\,y\leq v$: 
$$
\frac{1}{2}\left[\left|\frac{x-3y}{5}-\frac{u-3v}{5}\right|+ \left|\frac{y-3x}{5}-\frac{v-3u}{5}\right|\right] \leq \frac{2}{5}\left[\left|x-u\right|+\left|y-v\right|\right]
$$
$$
<\frac{2}{5} \cdot 2(\epsilon+\delta(\epsilon)) <\epsilon, 
$$
which holds if we simply take $\delta(\epsilon)<\frac{5}{4} \epsilon$. Thus, condition \eqref{Samet1} holds. Note also that $x_0=-3,\,y_0=3$ satisfy \eqref{mic}.

So Theorem \ref{th2} can be applied to $F$ in this example to conclude that $F$ has a (unique) coupled fixed point $(0,0)$, while Theorem \ref{th1} cannot be applied since \eqref{Samet} is not satisfied.
\end{example} 

\begin{remark} \em
One can prove that the coupled fixed point ensured by Theorem \ref{th2} is in fact unique, like in Example \ref{ex1}, provided that:
\end{remark}
\textnormal{every pair of elements in} $X^2$ \textnormal{has either a lower bound or an upper bound},
which is known, see \cite{Bha}, to be equivalent to the following condition: for all $Y=(x,y),\,\overline{Y}=(\overline{x},\overline{y})\in X^2$,
\begin{equation} \label{eq-7}
\exists Z=(z_1,z_2)\in X^2 \textnormal{ that is comparable to }  Y \textnormal{ and } \overline{Y}.
\end{equation}

\begin{theorem} \label{th3}
Adding  condition  \eqref{eq-7}  to  the  hypotheses  of   Theorem \ref{th2},  we  obtain  the
uniqueness of the coupled fixed point of  $F$.
\end{theorem}

\begin{proof} By Theorem \ref{th2} there exists a coupled fixed point $(\overline{x},\overline{y})$. 
In search for a contradiction, assume that $Z^*=(x^*,y^*)\in X^2$ is a coupled fixed point of $F$, different from $\overline{Z}=(\overline{x},\overline{y})$. This means that $d_2(Z^*,\overline{Z})>0.$ We discuss two cases:

Case 1. $Z^*$ is comparable to $\overline{Z}$.

As $Z^*$ is comparable to $\overline{Z}$ with respect to the ordering in $X^2$, by taking in \eqref{contr} $Y=Z^*$ and $V=\overline{Z}$ (or $V=Z^*$ and $Y=\overline{Z}$), we obtain
$$
d_2(Z^*,\overline{Z})=d_2(T(Z^*),T(\overline{Z}))<d_2(Z^*,\overline{Z}),
$$
a contradiction.

Case 2. $Z^*$ and $\overline{Z}$ are not comparable.

In this case, there exists an upper bound or a lower bound $Z=(z_1,z_2)\in X^2$ of $Z^*$ and $\overline{Z}$. Then, in view of the monotonicity of $T$, $T^n(Z)$ is comparable to $T^n(Z^*)=Z^*$ and to $T^n(\overline{Z})=\overline{Z}$. Assume, without any loss of generality, that $\overline{x}<z_1,\,\overline{y}\geq z_2$ and $x^*<z_1,\,y^*\geq z_2$, which means
$\overline{Z}<Z$ and $Z^*<Z$. By the monotonicity of $T$, we have
$$
T^n(\overline{Z})\leq T^n(Z) \textnormal{ and } T^n(Z^*)\leq T^n(Z), \,\textnormal{ for all } n\geq 0
$$
Note that, like in the proof of Theorem \ref{th2}, condition \eqref{contr} implies the strict contractive condition
\begin{equation} \label{eq-5.6}
d_2(T(Y),T(Z))< d_2 (Y, Z),\, Y>Z,
\end{equation}
Take now $Y=\overline{Z} > Z=V$ in \eqref{eq-5.6} to obtain
\begin{equation} \label{eq-5.7}
d_2(T^{n+1}(\overline{Z}),T^{n+1}(Z))< d_2(T^{n}(\overline{Z}),T^{n}(Z)),\,n\geq 1,
\end{equation}
which shows that the sequence of nonnegative numbers $\{\eta_n\}_{n=0}^{\infty}$ given by
\begin{equation} \label{eq-5.8}
\eta_n=d_2(T^{n}(\overline{Z}),T^{n}(Z)),\,n\geq 1,
\end{equation}
is non-increasing, hence convergent to some $\epsilon\geq 0$.  

We now prove that necessarily $\epsilon=0$. Suppose, to the contrary, that $\epsilon>0$. Then, there exists a positive integer $p$ such that
$$
\epsilon\leq \eta_p<\epsilon+\delta(\epsilon),
$$
which, by the Meir-Keeler condition \eqref{contr}, yields
$$
\eta_{p+1}=d_2(T^{p+1}(\overline{Z}),T^{p+1}(Z))<\epsilon,
$$
a contradiction, since $\{\eta_n\}_{n=0}^{\infty}$  converges non-increasingly to $\epsilon$. Therefore $\epsilon=0$, that is, 
\begin{equation} \label{eq-5.9}
\lim_{n\rightarrow \infty} d_2(T^{n}(\overline{Z}),T^{n}(Z))=0.
\end{equation}
Similarly, one obtains
\begin{equation} \label{eq-5.10}
\lim_{n\rightarrow \infty} d_2(T^{n}(Z^*),T^{n}(Z))=0.
\end{equation}
Now, by  \eqref{eq-5.9} and \eqref{eq-5.10}, we have
$$
d_2(Z^*,\overline{Z})=d_2(T^n(Z^*),T^n(\overline{Z}))\leq 
$$
$$
d_2(T^n(Z^*),T^n(Z))+d_2(T^n(Z),T^n(\overline{Z}))\rightarrow 0
$$
as $n\rightarrow \infty$,

\noindent which leads to the contradiction $0<d_2(Z^*,\overline{Z})\leq 0$.
\end{proof}

Similarly to \cite{Bha} and \cite{Sam}, by assuming a similar  condition to \eqref{eq-7}, but this time with respect to the ordered set $X$, that is, by assuming that every pair of elements of $X$ have either an upper bound or a lower bound in $X$, one can show that even the components of the coupled fixed points are equal. 
 
\begin{theorem} \label{th4}
In addition to the hypotheses of  Theorem \ref{th3}, suppose that every pair of elements of  $X$ has an upper bound or a lower bound in $X$. Then for the coupled fixed point $(\overline{x},\overline{y})$ we have $\overline{x} = \overline{y}$, that is, $F$ has a fixed point:
$$
F(\overline{x},\overline{x})=\overline{x}.
$$  
\end{theorem}

\begin{proof}
Let $(\overline{x},\overline{y})$ be a coupled fixed point of $F$ (ensured by Theorem \ref{th2}). Suppose, to the contrary, that $\overline{x} \neq \overline{y}$. Without any loss of generality, we can assume $\overline{x} > \overline{y}$. Then
$$
\overline{x}=F(\overline{x},\overline{y})>\overline{y}=F(\overline{y},\overline{x}).
$$
We consider again two cases.

Case 1. If $\overline{x},\overline{y}$ are comparable, then $
F(\overline{x},\overline{y})=\overline{x}$  is comparable to $\overline{y}=F(\overline{y},\overline{x})$ and hence, by taking $x:=\overline{x},\,y:=\overline{y},\,u:=\overline{y},\,v:=\overline{x},\,$ in \eqref{eq-5.6} one obtains
\begin{equation} \label{eq-9}
0<d(\overline{x},\overline{y})=d(F(\overline{x},\overline{y}),F(\overline{y},\overline{x}))< d(\overline{x},\overline{y}),
\end{equation}
a contradiction.

Case 2. If $\overline{x},\overline{y}$ are not comparable, then there exists a $z\in X$ comparable to $\overline{x}$ and $\overline{y}$. Suppose $\overline{x}\leq z$ and $\overline{y}\leq z$ (the other case is similar). Then in view of the order on $X^2$, it follows that 
$$
(\overline{x},\overline{y})\geq (\overline{x},z);\,(\overline{x},z)\leq (z,\overline{x});\,(z,\overline{x})\geq (\overline{y},\overline{x}),
$$
that is $(\overline{x},\overline{y}),\, (\overline{x},z);\,(\overline{x},z),\,(z,\overline{x});\,(z,\overline{x}),\, (\overline{y},\overline{x})$ are comparable in $X^2$. Now, similarly to the proof of Theorem \ref{th3}, we obtain that,  for any two comparable elements $Y,V$ in $X^2$, one has
\begin{equation} \label{eq-8}
\lim_{n\rightarrow \infty} d_2(T^{n}(Y),T^n(Z))=0.
\end{equation}
where $T$ was defined in the proof of Theorem \ref{th3}.

Now use \eqref{eq-8} for the comparable pairs $Y=(\overline{x},\overline{y}),\, V=(\overline{x},z);\,Y=(\overline{x},z),\,V=(z,\overline{x});\,Y=(z,\overline{x}),\, V=(\overline{y},\overline{x})$, respectively, to get

\begin{equation} \label{eq-10}
\lim_{n\rightarrow \infty} d_2(T^n(\overline{x},\overline{y}),T^n(\overline{x},z))=0,
\end{equation}

\begin{equation} \label{eq-11}
\lim_{n\rightarrow \infty} d_2(T^n(\overline{x},z),T^n(z,\overline{x})=0,
\end{equation}

\begin{equation} \label{eq-12}
\lim_{n\rightarrow \infty} d_2(T^n(z,\overline{x}),T^n(\overline{y},\overline{x}))=0.
\end{equation}
Now, by using the triangle inequality and \eqref{eq-10}, \eqref{eq-11}, \eqref{eq-12}, one has
$$
d(\overline{x},\overline{y})=\frac{d(\overline{x},\overline{y})+d(\overline{x},\overline{y})}{2}=d_2((\overline{x},\overline{y}),(\overline{y},\overline{x}) )=d_2(T^n(\overline{x},\overline{y}),T^n(\overline{y},\overline{x}))
$$
$$
\leq d_2(T^n(\overline{x},\overline{y}),T^n(\overline{x},z))+d_2(T^n(\overline{x},z),T^n(z,\overline{x}))+
$$
$$
+d_2(T^n(z,\overline{x}),T^n(\overline{y},\overline{x})) \rightarrow 0 \,\textnormal{ as } n\rightarrow\infty,
$$
which shows that $d(\overline{x},\overline{y})=0$, that is $\overline{x}=\overline{y}$.
\end{proof}

Similarly, one can obtain the same conclusion under the following alternative assumption.

\begin{theorem} \label{th6}
In addition to the e of  Theorem \ref{th3}, suppose that $x_0,y_0 \in X$ are comparable. Then for the coupled fixed point $(\overline{x},\overline{y})$ we have $\overline{x} = \overline{y}$, that is, $F$ has a fixed point:
$$
F(\overline{x},\overline{x})=\overline{x}.
$$   
\end{theorem}

\begin{remark}
\em Note that our contractive condition \eqref{Samet1} is symmetric, while the contractive condition \eqref{Samet} used in \cite{Sam} is not. Our generalization is based in fact on the idea of making the last one symmetric, which is very natural, as the great majority of contractive conditions in metrical fixed point theory are symmetric, see \cite{Ber07} and \cite{Rus2}.
\end{remark}

\begin{remark}
\em Note also that if $F$ satisfies the contractive condition in \cite{Bha}, that is, there exists a constant $k \in \left[0,1\right)$  with                                                                                   
$$
	d\left(F\left(x,y\right),F\left(u,v\right)\right) \leq \frac{k}{2}\left[d\left(x,u\right) + d\left(y,v\right)\right],\textnormal{ for each }x \geq u, y \leq  v.
$$
then, as pointed out by Proposition 2.1 in \cite{Sam}, $F$ also satisfies the contractive condition \eqref{Samet} and hence \eqref{Samet1}.

This follows by simply taking $\delta(\epsilon)=\left(\dfrac{1}{k}-1\right)\epsilon$.

In view of the results in \cite{Lim} and \cite{Suz}, the coupled fixed point theorems established in the present paper are also generalizations of all results in  \cite{Aga}, \cite{Ran}, \cite{Meir} and \cite{Sam}.
\end{remark}

\vskip 0.5 cm {\it 

Department of Mathematics and Computer Science

North University of Baia Mare

Victoriei 76, 430122 Baia Mare ROMANIA

E-mail: vberinde@ubm.ro}
\end{document}